\documentstyle{amsppt}

\topmatter
\title
ON THE CYCLIC SUBGROUP SEPARABILITY OF FREE PRODUCTS OF TWO GROUPS WITH AMALGAMATED SUBGROUP
\endtitle
\author
E. V. Sokolov
\endauthor
\affil
Ivanovo State University \\
Ermak st., 39, 153025, Ivanovo, Russia \\
{\it E-mail:} ev-sokolov\@yandex.ru
\endaffil
\address
Ivanovo State University
Ermak st., 39, 153025, Ivanovo, Russia
\endaddress
\email
ev-sokolov\@yandex.ru
\endemail
\abstract
Let $G$ be a free product of two groups with
amalgamated subgroup, $\pi $ be either the set of all prime
numbers or the one-element set \{$p$\} for some prime number $p$.
Denote by $\Sigma $ the family of all cyclic subgroups of group
$G$, which are separable in the class of all finite $\pi $@-groups.
\endgraf
Obviously, cyclic subgroups of the free factors, which aren't
separable in these factors by the family of all normal subgroups
of finite $\pi $@-index of group $G$, the subgroups conjugated
with them and all subgroups, which aren't $\pi ^{\prime}
$@-isolated, don't belong to $\Sigma $. Some sufficient
conditions are obtained for $\Sigma $ to coincide with the family
of all other $\pi ^{\prime} $@-isolated cyclic subgroups of group~$G$.
\endgraf
It is proved, in particular, that the residual $p$@-finiteness of
a free product with cyclic amalgamation implies the
$p$@-separability of all $p^{\prime} $@-isolated cyclic subgroups
if the free factors are free or finitely generated residually
$p$@-finite nilpotent groups.
\endabstract
\endtopmatter

\rightheadtext{ON THE CYCLIC SUBGROUP SEPARABILITY}

{\bf 1.~Introduction. Main results}

\bigskip

Let ${\Cal K}$ be a class of groups. We recall (see [1]), that a
subgroup $F$ of a group $G$ is said to be separable by the groups
of class ${\Cal K}$ if, to any element $g \in G\smallsetminus F$,
there exists a homomorphism $\psi $ of group $G$ onto a group of
${\Cal K}$ such that $g\psi \notin F\psi $. Group $G$ is called
residually ${\Cal K}$ if it's trivial subgroup is separable by the
groups of class ${\Cal K}$. If class ${\Cal K}$ coincides with the
class of all finite groups, then we shall say about residual
finiteness and about finite separability of subgroups. A group,
all cyclic subgroups of which are finitely separable, is called
$\pi _{c}$@-group.

If $\Psi $ is a family of normal subgroups of group $G$, then we
shall say also that subgroup $F$ of group $G$ is separable by the
subgroups of $\Psi $ if $\bigcap _{N \in \Psi }FN = F$. Thus, the
separability of subgroup $F$ in class ${\Cal K}$ is equivalent to
the separability of $F$ by the family of all normal subgroups of
group $G$, the factor-groups by which belong to ${\Cal K}$.

It is obvious that the separability of all cyclic subgroups of
group $G$ in class ${\Cal K}$ implies the ``residually ${\Cal
K}$'' property of $G$. The converse, in general, isn't true, and
so the problem arises to describe all cyclic subgroups of a
residually ${\Cal K}$ group, which are separable by the groups of
${\Cal K}$. The case considered in the present paper is that
${\Cal K}$ coincides with the class of all finite groups or all
finite $p$@-groups and group $G$ is a free product of two groups
with amalgamated subgroup.

Let $A$ and $B$ be some groups, $H$ be a subgroup of $A$, $K$ be a
subgroup of $B$ and let $\varphi $:~$H \to K$ be an isomorphism.
Let $G = (A * B; H = K, \varphi )$ be the free product of groups
$A$ and $B$ with subgroups $H$ and $K$ amalgamated according to
isomorphism $\varphi $. Obviously, an arbitrary cyclic subgroup of
group $G$ conjugated with a subgroup of one of the free factors
$A$ and $B$, that isn't separable in class ${\Cal K}$, will not be
separable by the groups of ${\Cal K}$ itself. Thus, our task is to
determine which of the remaining cyclic subgroups are separable in
${\Cal K}$ and, in particular, to find the conditions for them all
being separable by the groups of ${\Cal K}$.

\medskip

More precisely, the problem may be formulated as follows. Let
$\Delta _{A}$ and $\Delta _{B}$ be the families of all cyclic
subgroups of groups $A$ and $B$, respectively, which aren't
finitely separable in these groups. It is necessary to find the
conditions guaranteeing the truth of the following statement:

\medskip

$(*)$ \quad {\it An arbitrary cyclic subgroup of group $G$, which isn't
conjugate with any subgroup of the family $\Delta _{A} \cup \Delta
_{B}$, is finitely separable.}

\medskip

We note that in a series of papers (see, e.~g., [2]), dealing with
generalized free products of two $\pi _{c}$@-groups, the special
case of this task was considered, when the family $\Delta _{A}
\cup \Delta _{B}$ was empty.

If, as above, $H$ and $K$ are subgroups of groups $A$ and $B$,
respectively, and $\varphi $:~$H \to K$ is an isomorphism, then,
following G.~Baumslag [3], we shall call subgroups $R \leqslant A$
and $S \leqslant B$ $(H, K, \varphi )$@-compatible if $(R \cap
H)\varphi = S \cap K$. Let denote by $\Omega $ the family of all
pairs of normal $(H, K, \varphi )$@-compatible subgroups of finite
index of groups $A$ and $B$ and by $\Omega _{A}$ and $\Omega _{B}$
it's projections onto groups $A$ and $B$.

It is easy to see that, if $N$ is an arbitrary normal subgroup of
finite index of group $G$, then the pair $(A \cap N, B \cap N)$
belongs to family $\Omega $. It follows from this remark that a
finitely separable cyclic subgroup of group $G$ contained in one
of groups $A$ and $B$ is separable by the subgroups of families
$\Omega _{A}$ or $\Omega _{B}$, respectively.

Let now $\Lambda _{A}$ and $\Lambda _{B}$ denote the families of
all cyclic subgroups of groups $A$ and $B$, which aren't separable
by the subgroups of families $\Omega _{A}$ or $\Omega _{B}$. Then
the condition just stated can be formulated in the form of the
following

\proclaim{\indent Proposition~1.1} If a cyclic subgroup of group
$G$ is finitely separable, then it conjugates with no subgroup of
the family $\Lambda _{A} \cup \Lambda _{B}$.
\endproclaim

We note now that inclusions $\Delta _{A} \subseteq \Lambda _{A}$
and $\Delta _{B} \subseteq \Lambda _{B}$ take place, since the
finite separability of a subgroup of a given group means precisely
the separability by the family of all normal subgroups of finite
index. Thus, statement~$(*)$ is equivalent to simultaneous
realizability of the next two ones:

\noindent
a)~$\Delta _{A} = \Lambda _{A}$ and $\Delta _{B} =
\Lambda _{B}$, and

\noindent
b)~an arbitrary cyclic subgroup of group $G$, that isn't
conjugate with any subgroup of $\Lambda _{A} \cup \Lambda _{B}$,
is finitely separable.

The following statement, the first of the main results of the
paper, gives a sufficient condition for the second claim to be
true.

\proclaim{\indent Theorem~1.2} Let family $\Omega _{A}$ be an
$H$@-filtration and family $\Omega _{B}$ be a $K$@-filtration.
Then an arbitrary cyclic subgroup of group $G$, that conjugates
with no subgroup of $\Lambda _{A} \cup \Lambda _{B}$, is finitely
separable.
\endproclaim

We recall (see [3]) that a family $\Psi $ of normal subgroups of a
group $X$ is said to be a $Y$@-filtration, where $Y$ is a subgroup
of $X$, if $\bigcap_{N \in \Psi }N = 1$ and $Y$ is separable by the
subgroups of $\Psi $. Proposition~2 of paper [3] asserts that, if
family $\Omega _{A}$ is an $H$@-filtration and family $\Omega
_{B}$ is a $K$@-filtration, then $G$ is a residually finite group.
Thus, theorem~1.2 may be considered as a generalization of this
statement.

Having slightly increased our restrictions, we may obtain the
maximal property $(*)$ for the family of finitely separable
cyclic subgroups of group $G$.

\proclaim{\indent Theorem~1.3} Let groups $A$ and $B$ be
residually finite, subgroups $H$ and $K$ be finitely separable in
the free factors and, to any two normal subgroups of finite index
$M \leqslant A$ and $N \leqslant B$ there exists a pair of
subgroups $(R, S) \in \Omega $ such that $R \leqslant M$ and $S
\leqslant N$. Then group $G$ satisfies condition~$(*)$. In
particular, if $A$ and $B$ are $\pi _{c}$@-groups, then $G$ is
also a $\pi _{c}$@-group.
\endproclaim

Indeed, any subgroup $F$ of group $A$ or group $B$, which is
finitely separable in $A$ or $B$, turns out separable by the
subgroups of $\Omega _{A}$ or $\Omega _{B}$, respectively, in this
case. Therefore, in particular, $\Delta _{A} = \Lambda _{A}$ and
$\Delta _{B} = \Lambda _{B}$. Besides, groups $A$ and $B$ being
residually finite, $\Omega _{A}$ is an $H$@-filtration and $\Omega
_{B}$ is a $K$@-filtration. The desired claim follows now from
theorem~1.2.

We note that theorems~1.2 and 1.3 are a generalization of the
results obtained by G.~Kim [2, theorem~1.1 and proposition~1.2]
for generalized free products of $\pi _{c}$@-groups.

\medskip

Let turn now to description of cyclic subgroups of group $G = (A
* B; H = K, \varphi )$, which are separable in the class of
finite $p$@-groups (or, briefly, $p$@-separable).

We remind, first of all, that a subgroup $Y$ of a group $X$ is
called $p^{\prime} $@-isolated if, for any element $g \in Y$ and
for any prime number $q$, which doesn't equal $p$, $g^{q} \in Y$
implies $g \in Y$. It is easy to see that every $p$@-separable
subgroup must be $p^{\prime} $@-isolated, and so the original task
takes the following form.

Let $\Delta _{A}^{p} $ and $\Delta _{B}^{p} $ be the families of all
$p^{\prime} $@-isolated cyclic subgroups of groups $A$ and $B$, respectively,
which aren't $p$@-separable in these groups. It is necessary to find the
conditions guaranteeing the truth of the following statement:

\medskip

$(**)$ \quad {\it An arbitrary $p^{\prime} $@-isolated cyclic subgroup of
group $G$, which isn't conjugate with any subgroup of $\Delta
_{A}^{p} \cup \Delta _{B}^{p} $, is $p$@-separable.}

\medskip

We remark that the $p$@-separability of all $p^{\prime} $@-isolated cyclic
subgroups of some group doesn't necessarily imply the residual
$p$@-finiteness of this group.

Let $\Omega ^{p}$ denotes the family of all ordered pairs $(A \cap
N, B \cap N)$, where $N$ is an arbitrary normal subgroup of group
$G$ of finite $p$@-index. Let also $\Omega _{A}^{p} $ and $\Omega
_{B}^{p} $ denote the families of the first and the second
components of elements of $\Omega ^{p}$. The next proposition is
obtained by E.~D.~Loginova in the paper [4].

\proclaim{\indent Proposition~1.4}~A pair of subgroups $(R, S)$
belongs to family $\Omega ^{p}$ if, and only if there exist
sequences of subgroups $R = R_{0} \leqslant \ldots \leqslant R_{m}
= A$, $S = S_{0} \leqslant \ldots \leqslant S_{n} = B$ such that:

1) $R_{i}$, $S_{j}$ are normal subgroups of groups $A$ and $B$,
respectively $(0 \leqslant i \leqslant m, 0 \leqslant j \leqslant
n)$;

2) $\vert R_{i + 1}/R_{i}\vert = \vert S_{j + 1}/S_{j}\vert = p$
$(0 \leqslant i \leqslant m-1, 0 \leqslant j \leqslant n-1)$;

3) isomorphism $\varphi $ maps the set \{$R_{i} \cap H$\} onto the set
\{$S_{j} \cap K$\}.
\endproclaim

Following to [4] we shall call subgroups $R$ and $S$ satisfying
the conditions of proposition~1.4 $(H, K, \varphi ,
p)$@-compatible.

Let us denote by $\Lambda _{A}^{p} $ and $\Lambda _{B}^{p} $ the
families of all $p^{\prime} $@-isolated cyclic subgroups of groups
$A$ and $B$, which aren't separable by the subgroups of $\Omega
_{A}^{p} $ and $\Omega _{B}^{p} $, respectively. Then, as above,
the inclusions $\Delta _{A}^{p} \subseteq \Lambda _{A}^{p} $,
$\Delta _{B}^{p} \subseteq \Lambda _{B}^{p} $ take place, and the
following proposition is true.

\proclaim{\indent Proposition~1.5}~If a $p^{\prime} $@-isolated
cyclic subgroup of group $G$ is $p$@-separable, then it conjugates
with no subgroup of the family $\Lambda _{A}^{p} \cup \Lambda
_{B}^{p} $.
\endproclaim

In the same paper [4] the analog of the mentioned above sufficient
condition by Baumslag is obtained: if family $\Omega _{A}^{p} $ is
an $H$@-filtration and family $\Omega _{B}^{p} $ is a
$K$@-filtration, then group $G$ is residually $p$@-finite. It
turns out that the statements, similar to theorems~1.2 and 1.3,
also take place.

\proclaim{\indent Theorem~1.6} Let family $\Omega _{A}^{p} $ be an
$H$@-filtration and family $\Omega _{B}^{p} $ be a
$K$@-filtration. Then an arbitrary $p^{\prime} $@-isolated cyclic
subgroup of group $G$, which conjugates with no subgroup of
$\Lambda _{A}^{p} \cup \Lambda _{B}^{p} $, is $p$@-separable.
\endproclaim

\proclaim{\indent Theorem~1.7}~Let groups $A$ and $B$ be
residually $p$@-finite, subgroups $H$ and $K$ be $p$@-separable in
the free factors and, to any two normal subgroups of finite
$p$@-index $M \leqslant A$ and $N \leqslant B$, there exists a
pair of subgroups $(R, S) \in \Omega ^{p}$ such that $R \leqslant
M$ and $S \leqslant N$. Then group $G$ satisfies condition~$(**)$.
\endproclaim

The last theorem is deduced from theorem~1.6 in exactly the same
way as theorem~1.3 from theorem~1.2.

\bigskip

{\bf 2.~Some applications}

\bigskip

Let $A$ be a free group with the set of free generators \{$a,
b$\}, $B$ be a free group with the set of free generators \{$c,
d$\}, and let $H$ be the subgroup of group $A$ generated by the
elements $a$ and $a_{1} = b^{ - 1}ab$, $K$ be the subgroup of
group $B$ generated by the elements $c$ and $c_{1} = d^{ -
1}c^{2}d$. It is obvious that the indicated generators of
subgroups $H$ and $K$ generate these subgroups freely, and so the
map, which associates $a$ with $c$ and $a_{1}$ with $c_{1}$,
defines an isomorphism $\varphi $ of subgroup $H$ onto subgroup
$K$.

Thus, the group $G = \langle a, b, c, d; a = c, b^{ - 1}ab = d^{
- 1}c^{2}d\rangle $ is a free product of groups $A$ and $B$ with
subgroups $H$ and $K$ amalgamated according to isomorphism
$\varphi $.

\proclaim{\indent Theorem~2.1}~An arbitrary cyclic subgroup of the
group $G = \langle a , b, c, d; a = c, b^{-1}ab =
d^{-1}c^{2}d\rangle $, that isn't conjugate with any subgroup of
$\Lambda _{A} \cup \Lambda _{B}$, is finitely separable. At the
same time families $\Lambda _{A}$ and $\Lambda _{B}$ aren't empty,
and family $\Omega _{B}$ isn't a $K$@-filtration.
\endproclaim

The given statement demonstrates that the sufficient condition
stated in theorem~1.2 isn't necessary. Besides, all finitely
generated subgroups of an arbitrary free group being finitely
separable, the first two conditions of theorem~1.3 are fulfilled
here. But $G$ isn't a $\pi _{c}$@-group. Thus, the third condition
of this theorem isn't true and hence doesn't follow, in general,
from the first two ones.

It isn't difficult to verify that the $(H, K, \varphi
)$@-compatibility of normal subgroups of finite $p$@-index implies
their $(H, K, \varphi , p)$@-compatibility in a free product of
two groups with cyclic amalgamation. It is easy to see also that,
to any element $g$ of a residually $p$@-finite group and to any
$p$@-number $x$, there exists a normal subgroup of finite
$p$@-index, which intersects with the cyclic subgroup generated by
$g$ at the subgroup $\langle g^{x}\rangle $.

Thus, if groups $A$ and $B$ are residually $p$@-finite, and
subgroups $H$ and $K$ are cyclic, then families $\Omega _{A}^{p} $
and $\Omega _{B}^{p} $ coincide with the families of all normal
subgroups of groups $A$ and $B$ of finite $p$@-index. This remark
results in the next statement following directly from theorem~1.7.

\proclaim{\indent Theorem~2.2} Let $A$ and $B$ be residually
$p$@-finite groups, $H$ and $K$ be cyclic subgroups, which are
$p$@-separable in the free factors. Then group $G$ is residually
$p$@-finite and satisfies condition~$(**)$.
\endproclaim

The same reasons are used in the proof of one more result.

\proclaim{\indent Theorem~2.3} Let $H$ and $K$ are infinite cyclic
subgroups, and their centralizers in groups $A$ and $B$,
respectively, don't contain elements of finite order. If group $G$
is residually $p$@-finite, then it satisfies condition~$(**)$.
\endproclaim

Let us formulate now two statements following directly from
theorems~2.2 and 2.3, respectively.

\proclaim{\indent Corollary~2.4}~Let $A$ and $B$ be finitely
generated residually $p$@-finite nilpotent groups (i.~e. their
torsion parts are $p$@-groups), $\langle h\rangle \leqslant A$
and $\langle k\rangle \leqslant B$ be maximal infinite cyclic
subgroups and $H = \langle h^{m}\rangle $, $K = \langle
k^{n}\rangle $. If $m$ and $n$ are $p$@-numbers, then all
$p^{\prime} $@-isolated cyclic subgroups of group $G$ are
$p$@-separable.
\endproclaim

\proclaim{\indent Corollary~2.5}~Let $A$ and $B$ be free groups,
$\langle h\rangle \leqslant A$ and $\langle k\rangle \leqslant B$
be maximal cyclic subgroups and $H = \langle h^{m}\rangle $, $K =
\langle k^{n}\rangle $. If $m = 1$ or $n = 1$ or $m$ and $n$ are
$p$@-numbers, then all $p^{\prime} $@-isolated cyclic subgroups of
group $G$ are $p$@-separable.
\endproclaim

It is proved in [5] and [6] that the conditions of corollaries~2.4
and 2.5 are necessary and sufficient for the residual
$p$@-finiteness of group $G$. So the only remark which is needed
for the proof is that all $p^{\prime} $@-isolated cyclic subgroups
of free and finitely generated nilpotent groups are $p$@-separable
(see [5] and [4], respectively).

We note that, as it followes from corollary~2.5, neither the
$p$@-separability of the amalgamated subgroups, nor even their
$p^{\prime} $@-isolation isn't the necessary condition for the
$p$@-separability of all $p^{\prime} $@-isolated cyclic subgroups
of group $G$.

The other applications of theorems~1.3 and 1.7 can be founded in
the author's papers [7] and [8].

\bigskip

{\bf 3.~The proof of theorems~1.2 and 1.6}

\bigskip

To any pair of subgroups $(R, S) \in \Omega $, the map $\varphi
_{R, S}: HR/R \to KS/S$, which associates an element $hR$, $h \in
H$, with the element $(h\varphi )S$, is correctly defined and
serves as an isomorphism of subgroups. Therefore we may construct
the group $G_{R, S} = (A/R * B/S; HR/R = KS/S, \varphi _{R, S})$.
The natural homomorphisms of group $A$ onto $A/R$ and of group $B$
onto $B/S$ are extendable to a homomorphism $\pi _{R, S}$ of group
$G$ onto group $G_{R, S}$.

It is well known that generalized free product of two finite
groups is residually finite and moreover a $\pi _{c}$@-group. So,
to any pair of subgroups $(R, S) \in \Omega $, $G_{R, S}$ is a
$\pi _{c}$@-group.

Generalized free product of two finite $p$@-groups isn't, in
general, residually $p$@-fini\-te. The corresponding criteria was
founded by G.~Higman in [9]. It follows directly from this
criteria and proposition~1.4 that, if $(R, S) \in \Omega $, then
the group $G_{R, S}$ is residually $p$@-finite if, and only if
$(R, S) \in \Omega ^{p}$. We'll show also that all $p^{\prime}
$@-isolated cyclic subgroups of group $G_{R, S}$ are
$p$@-separable in this case.

\proclaim{\indent Proposition~3.1} Let $A$ and $B$ be finite
groups. If $G$ is a residually $p$@-finite group, then all it's
$p^{\prime} $@-isolated cyclic subgroups are $p$@-separable.
\endproclaim

{\it Proof.}~For group $G$ is residually $p$@-finite, there exists
it's homomorphism onto a finite $p$@-group, the kernel of which
intersects trivially with the free factors and, because of known
theorem by H.~Neumann [10], is a free group. As it was noted
above, all $p^{\prime} $@-isolated cyclic subgroups of free group
are $p$@-separable, and so the desired claim results from the
following statement.

\proclaim{\indent Proposition~3.2} Let a group $X$ be an extension
of a group $Y$ by a finite $p$@-group and let all $p^{\prime}
$@-isolated cyclic subgroups of group $Y$ be $p$@-separable. Then
all $p^{\prime} $@-isolated cyclic subgroups of group $X$ are
$p$@-separable too.
\endproclaim

{\it Proof.}~Let $F$ be a $p^{\prime} $@-isolated cyclic subgroup
of group $X$, $g \in X\smallsetminus F$. It is sufficient for
proving to point out a normal subgroup $N$ of finite $p$@-index
such that $g \notin FN$.

If $g \notin FY$, then subgroup $Y$ is desired. So $g$ will be
considered to be an element of $FY$.

We write $g$ in the form $g = fy$, where $f \in F$, $y \in Y$.
Since $g \notin F$, $y \notin F \cap Y$.

Obviously, $F \cap Y$ is a $p^{\prime} $@-isolated cyclic subgroup
of group $Y$. Hence it is $p$@-separable in $Y$ and there exists a
normal subgroup $M$ of group $Y$ of finite $p$@-index such that $y
\notin (F \cap Y)M$. To every element $y \in Y$, the subgroup $y^{
- 1}My$ is included in $Y$, is normal and has finite $p$@-index in
this group. Owing to finiteness of the index $[X:Y]$, the number
of different subgroups of such form is also finite. Thus, their
intersection $N$, say, is a subgroup of finite $p$@-index of group
$Y$, normal in $X$.

If $g \in FN$, then $g = f^{\prime} u$ for some elements
$f^{\prime} \in F$, $u \in N$ and $f^{- 1}f^{\prime} = yu^{- 1}
\in F \cap Y$. But $y = (f^{- 1}f^{\prime} )u \in (F \cap Y)N
\subseteq (F \cap Y)M$ in this case, what contradicts the choice
of subgroup $M$. Thus, $g \notin FN$, and subgroup $N$ is
required.

\proclaim{\indent Proposition~3.3} Let $X$ be a residually
$p$@-finite group, $g \in X$ be an element of infinite order. The
subgroup $\langle g\rangle $ isn't $p^{\prime} $@-isolated if, and
only if there exist an element $h \in X$ and a prime number $q$,
which doesn't equal $p$, such that $g = h^{q}$.
\endproclaim

{\it Proof.}~The sufficiency of this condition is obvious, we'll
show it's necessity.

Let $f \in X\smallsetminus \langle g\rangle $ be such an element
that $f^{q} \in \langle g\rangle $ for some prime number $q$,
which isn't equal to $p$, and $f^{q} = g^{k}$.

The residual $p$@-finiteness of group $X$ results that the
centralizer $C(g)$ of element $g$ in group $X$ is a $p$@-separable
subgroup and therefore a $p^{\prime} $@-isolated one. Hence $f \in
C(g)$.

If we suppose that $k = qk^{\prime} $, then $1 = f^{q}g^{-
qk^{\prime} } = (fg^{- k^{\prime}})^{q}$, and, owing to the
residual $p$@-finiteness of $X$, $f = g^{k^{\prime}} $. We obtain
a contradiction with the choice of element $f$. Thus, $(k, q) = 1$
and $ku + qv = 1$ for some integer numbers $u$ and $v$. From this
it follows that $g = g^{ku + qv} = f^{qu}g^{qv} =
(f^{u}g^{v})^{q}$, as claimed.

\medskip

We shall carry out {\bf the proof of theorems~1.2 and 1.6}
simultaneously and say about separability and compatibility of
subgroups without specifying of the concrete class of groups.

Let $h$ and $g$ be arbitrary elements of group $G$ such that $g
\ne 1$, the cyclic subgroup $\langle g\rangle $ conjugates with no
subgroup of $\Lambda _{A} \cup \Lambda _{B}$ (respectively, is
$p^{\prime} $@-isolated and conjugates with no subgroup of
$\Lambda _{A}^{p} \cup \Lambda _{B}^{p} $), and $h \notin \langle
g\rangle $. Let also $h = h_{1}h_{2} \ldots h_{m}$, $g =
g_{1}g_{2} \ldots g_{n}$ be reduced forms of elements $h$ and $g$.
Applying an appropriate inner automorphism of group $G$ we may
consider element $g$ as cyclically reduced.

To find a homomorphism $\theta $ of group $G$ onto a finite group
mapping $h$ to an element, which doesn't belong to $\langle
g\theta \rangle $, it is sufficient to point out a pair of
subgroups $(R, S) \in \Omega $ satisfying the property $h\pi
_{R, S} \notin \langle g\pi _{R, S}\rangle $. Since $G_{R, S}$ is
a $\pi _{c}$@-group, the homomorphism $\pi _{R, S}$ can be
extended to the desired homomorphism $\theta $.

We may use the same idea for constructing a homomorphism of group
$G$ onto a finite $p$@-group, but with minor restriction.

Indeed, if subgroups $R$ and $S$ are compatible, then all cyclic
subgroups of the free factors of the group $G_{R, S}$ (which are
finite $p$@-groups) are $p^{\prime} $@-isolated. Therefore, if $n
= 1$, i.~e. $g \in A \cup B$, and if we succeed to point out such
a pair of subgroups $(R, S) \in \Omega ^{p}$ that $h\pi _{R, S}
\notin \langle g\pi _{R, S}\rangle $, then the existence of the
required homomorphism follows from proposition~3.1.

But if $n \geqslant 2$, then the mere presence of a pair of
compatible subgroups $R$ and $S$ satisfying the property $h\pi
_{R, S} \notin \langle g\pi _{R, S}\rangle $ may turn out
insufficient, because the subgroup $\langle g\pi _{R, S}\rangle $
need not be $p^{\prime} $@-isolated in the group $G_{R, S}$ (the
corresponding example is given at the end of the proof). To make
use of proposition~3.1 in this case we shall find such a pair of
subgroups $(R, S) \in \Omega ^{p}$ that the image of $h$ under the
action of homomorphism $\pi _{R, S}$ doesn't belong to some
$p^{\prime} $@-isolated cyclic subgroup including the subgroup
$\langle g\pi _{R, S}\rangle $.

Let, at first, $n = 1$, and let $g \in A$ for definiteness.

By the condition the subgroup $\langle g\rangle $ is separable by
the subgroups of family $\Omega _{A}$ (respectively, of family
$\Omega _{A}^{p} $). Therefore, if $h \in A$, there exists a pair
of compatible subgroups $R$ and $S$ such that $h \notin \langle
g\rangle R$, and hence $h\pi _{R, S} \notin \langle g\pi
_{R, S}\rangle $.

Let $h \notin A$. Then $h \in B\smallsetminus K$ if $m = 1$ or
every syllable $h_{i}$ of it's reduced form belongs to one of the
free factors but isn't contained in the amalgamated subgroup if
$m>1$. So, to every $i$ ($1 \leqslant i \leqslant m$), we can
point out a pair of compatible subgroups $R$ and $S$ such that
$h_{i} \notin HR_{i}$ if $h_{i} \in A$ and $h_{i} \notin KS_{i}$
if $h_{i} \in B$. Let $R = \bigcap R_{i}$, $S = \bigcap S_{i}$.

It is easy to see that subgroups $R$ and $S$ are compatible,
$l(h\pi _{R, S}) = l(h)$ (here $l( \cdot )$ denotes syllable
length), and, if $m = 1$, then $h\pi _{R, S} \in B\pi
_{R, S}\smallsetminus K\pi _{R, S}$. Thus, $h\pi _{R, S} \notin
\langle g\pi _{R, S}\rangle $ in this case too.

Let now $n \geqslant 2$. We find, as above, a pair of compatible
subgroups $R$ and $S$ such that $l(h\pi _{R, S}) = l(h)$ and
$l(g\pi _{R, S}) = l(g)$. Obviously, the form of the element $g\pi
_{R, S}$ is cyclically reduced as before.

We shall finish the proof of theorem~1.2 at first.

It is not difficult to show that, for any two elements $u, v \in
G_{R, S}$, if one of these elements is cyclically reduced and $v
\in \langle u\rangle $, then the other element is also cyclically
reduced and $l(u)\vert l(v)$. Hence, if $n$ doesn't divide $m$,
then $h\pi _{R, S} \notin \langle g\pi _{R, S}\rangle $.

Let $m = nk$ for some positive $k$. Since $h \notin \langle
g\rangle $, then $h \ne g^{\pm k}$ and, because of the residual
finiteness of group $G$, there exists it's normal subgroup $L$ of
finite index not containing the elements $h^{ - 1}g^{k}$ and $h^{
- 1}g^{ - k}$. Putting $R^{\prime} = R \cap L$, $S^{\prime} = S
\cap L$ we have $l(h\pi _{R^{\prime}, S^{\prime}} ) = l(h)$,
$l(g\pi _{R^{\prime} , S^{\prime}} ) = l(g)$, and $h\pi
_{R^{\prime} , S^{\prime} } \ne (g\pi _{R^{\prime} , S^{\prime}}
)^{\pm k}$, whence follows that $h\pi _{R^{\prime} , S^{\prime}}
\notin \langle g\pi _{R^{\prime} , S^{\prime}} \rangle $.

Thereby, theorem~1.2 is proved, and we turn to the proof of
theorem~1.6.

Obviously, $l(g\pi _{R^{\prime} , S^{\prime}} ) = l(g)>1$ for any
pair of compatible subgroups $R^{\prime} $ and $S^{\prime} $,
which are included in $R$ and $S$, respectively. Applying
proposition~3.3 it isn't difficult to see that the subgroup
$\langle g\pi _{R^{\prime} , S^{\prime}} \rangle $ is contained in
some $p^{\prime} $@-isolated cyclic subgroup $F_{R^{\prime}
, S^{\prime}} $, it's index in this subgroup being mutually
distinct with $p$. We shall prove that subgroups $R^{\prime} $ and
$S^{\prime} $ can be chosen in a such way that the element $h\pi
_{R^{\prime} , S^{\prime}} $ doesn't belong to $F_{R^{\prime}
, S^{\prime}} $.

Let write the number $n$ in the form $n = p^{l}n^{\prime} $, where
($n^{\prime} , p) = 1$, and consider the two cases.

Case~1. $n$ doesn't divide $mn^{\prime} $.

Suppose that $h\pi _{R, S} \in F_{R, S}$. It is clear that the
index of the subgroup $\langle g\pi _{R, S}\rangle $ in group
$F_{R, S}$ divides $n^{\prime} $, and so $(h\pi
_{R, S})^{n^{\prime}} \in \langle g\pi _{R, S}\rangle $. But this
contradicts the supposition that $n$ doesn't divide $mn^{\prime}$.
Thus, $h\pi _{R, S} \notin F_{R, S}$.

Case~2. $mn^{\prime} = nk$ for some positive $k$.

Since the subgroup $\langle g\rangle $ is $p^{\prime} $@-isolated
in $G$ and $h \notin \langle g\rangle $, then $h^{n^{\prime}} \ne
g^{\pm k}$. The residual $p$@-finiteness of group $G$ results that
there exists a normal subgroup $L$ of group $G$ of finite
$p$@-index such that $h^{-n^{\prime}} g^{k}$, $h^{-n^{\prime}}
g^{-k} \notin L$. Let $R^{\prime} = R \cap L$, $S^{\prime} = S
\cap L$.

Then $(h\pi _{R^{\prime} , S^{\prime}} )^{n^{\prime}} \ne (g\pi
_{R^{\prime} , S^{\prime}} )^{\pm k}$, and so $(h\pi _{R^{\prime}
, S^{\prime}} )^{n^{\prime}} \notin \langle g\pi _{R^{\prime}
, S^{\prime}} \rangle $. It follows, as above, that $h\pi
_{R^{\prime} , S^{\prime}} \notin F_{R^{\prime} , S^{\prime}} $,
and the proof is finished.

\medskip

Let us make a remark now in connection with the given proof. Let
$F$ be a cyclic subgroup of group $G$ generated by a cyclically
reduced element $g$ of a syllable length greater than 1. It is
interesting that, even if subgroup $F$ is $p$@-separable, it may
be impossible to find such a pair of subgroups $(R, S) \in \Omega
^{p}$ that the element $g\pi _{R, S}$ has a reduced form of a
non-unit length, as before, and at the same time the subgroup
$\langle g\pi _{R, S}\rangle $ is $p^{\prime} $@-isolated in
$G_{R, S}$.

Let $G = \langle a, b; a^{p} = b^{p}\rangle $ and $g =
(ab)^{q}a^{p}$, where $p$,~$q$ are different prime numbers. It is
easy to see that the subgroup $\langle g\rangle $ is $p^{\prime}
$@-isolated in $G$ and hence is $p$@-separable in $G$ by virtue of
theorem~2.2.

From the other hand, for every pair of subgroups $(R, S) \in
\Omega ^{p}$, where $R \ne A$ and $S \ne B$, the group $G_{R, S}$
has the presentation $\langle a, b; a^{p^{n}} = b^{p^{n}} = 1,
a^{p} = b^{p}\rangle $ for a convenient natural $n$. Let $h =
aba^{px_{n}} $, where $x_{n}$ is a solution of the congruence $qx
\equiv 1~(mod~p^{n})$. Then, obviously, $h\pi _{R, S} \notin
\langle g\pi _{R, S}\rangle $ while $(h\pi _{R, S})^{q} \in
\langle g\pi _{R, S}\rangle $, and, thus, the subgroup $\langle
g\pi _{R, S}\rangle $ isn't $p^{\prime} $@-isolated in $G_{R, S}$.

\bigskip

{\bf 4.~The proof of theorems~2.1 and 2.3}

\bigskip

{\bf The proof of theorem~2.1.} We put $t = bd^{- 1}$ and then use
the obvious Tietze transformations to convert the presentation $G
= \langle a, b, c, d; a = c, b^{- 1}ab = d^{ - 1}c^{2}d\rangle $
of group $G$ to the presentation $G = \langle a, b, t; t^{- 1}at =
a^{2}\rangle $, which means that group $G$ is the ordinary free
product of the group $C = \langle a, t; t^{- 1}at = a^{2}\rangle $
and an infinite cyclic subgroup with generator $b$.

Owing to the residual finiteness of group $G$ and theorem~1.3 a
cyclic subgroup of group $G$ isn't finitely separable if, and only
if it conjugates with a subgroup of $\Delta _{C}$. It is well
known that family $\Delta _{C}$ consists of those subgroups of
group $C$, which conjugate with the subgroups generated by the
elements of form $a^{k}$.

We shall prove now the two auxiliary statements.

\proclaim{\indent Proposition~4.1} If a normal subgroup $M$ of
finite index of group $A$ (of group $B$) belongs to family $\Omega
_{A}$ (respectively, to family $\Omega _{B}$), then the order of
element $a$ (respectively, of element $c$) modulo subgroup $M$ is
an odd number.
\endproclaim

{\it Proof.} Let a normal subgroup $M$ of finite index of group
$A$ and a normal subgroup $N$ of finite index of group $B$ are
$(H, K, \varphi )$@-compatible. We put $H \cap M = U$ and $K \cap
N = V$, so that $U\varphi = V$.

It is obvious that the orders of elements $a$ and $a_{1}$ modulo
subgroup $M$ must coincide, and, the factor-group $H/U$ being
embeddable naturally to the factor-group $A/M$, the orders of
these elements modulo subgroup $U$ coincide too.

Considering the images according to isomorphism $\varphi $ we get
coincidence of the orders of elements $c$ and $c_{1}$ of group $K$
modulo subgroup $V$, and so coincidence of the orders of these
elements modulo subgroup $N$. It follows that elements $c$ and
$c^{2}$ have the same order, and therefore the order of element
$c$ modulo subgroup $N$ is an odd number.

Since the element $cN$ corresponds to the element $aM$ under the
isomorphism of the subgroup $HM/M$ of $A/M$ onto the subgroup
$KN/N$ of $B/N$, which is induced by isomorphism $\varphi $, the
order of element $a$ modulo subgroup $M$ is an odd number too.

\proclaim{\indent Proposition~4.2} A cyclic subgroup of group $A$
(of group $B$) belongs to family $\Lambda _{A}$ (respectively, to
family $\Lambda _{B}$) if, and only if it conjugates with a
subgroup generated by an element $a^{2k}$ (respectively, $c^{2k}$)
for some $k \ne 0$.
\endproclaim

{\it Proof.} The cyclic subgroup $F$ of group $A$ generated by the
element $a^{2k}$ doesn't contain the element $a^{k}$. Let $M$ be
an arbitrary subgroup of family $\Omega _{A}$. Owing to
proposition~4.1 the order $m$ of element $a$ modulo subgroup $M$
is an odd number, and so the congruence $2l \equiv 1~(mod~m)$ is
solvable for some integer number $l$. Therefore $a \equiv a^{2l}~
(mod~M)$, $a^{k} \equiv (a^{2k})^{l}~(mod~M)$, whence $a^{k} \in
FM$. Thus, subgroup $F$ isn't separable by family $\Omega _{A}$.

Conversely, the elements $a$ and $a^{2}$ conjugated in $G$, an
arbitrary cyclic subgroup $F$, which is contained in $A$ and
conjugates with no subgroup generated by an element $a^{2k}$, is
finitely separable in $G$, and, in accordance with
proposition~1.1, is separable by family $\Omega _{A}$.

The argument for group $B$ is analogous.

\medskip

Since the elements $a$ and $a^{2}$ are conjugated in group $G$,
the statement of theorem follows directly from propositions~4.1
and 4.2.

\medskip

{\bf The proof of theorem~2.3.} If families $\Omega _{A}^{p} $ and
$\Omega _{B}^{p} $ are an $H$@- and a $K$@-filtration, respectively,
the desired claim results from theorem~1.6. So this condition will
be considered to be false.

Let, for definiteness, family $\Omega _{A}^{p} $ be not an
$H$@-filtration. Owing to the residual $p$@-finiteness of group
$G$, this means that subgroup $H$ isn't separable by the subgroups
of family $\Omega _{A}^{p} $, and hence there exists an element $f
\in A\smallsetminus H$ moving to $H$ under the action of any
homomorphism of group $G$ onto a finite $p$@-group (we will denote
the family of all such homomorphisms by $\Psi $). Let us remark
that family $\Omega _{B}^{p} $ must be a $K$@-filtration then:
otherwise there exists an element $g$ of the set $B\smallsetminus
K$ with the analogous property, and the commutator $[f, g]$ turns
out a non-trivial element of group $G$, mapped to unit under any
homomorphism $\psi \in \Psi $.

Let further $h$ and $k$ be a generators of subgroups $H$ and $K$,
respectively, and $h\varphi = k$. First of all we'll show that,
to any natural $p$@-number $n$, there exists such an element
$f_{n} \in A\smallsetminus H$ that $f_{n}\psi \in H^{n}\psi $
under every homomorphism $\psi \in \Psi $.

Let $f \in A\smallsetminus H$ be an element moving to $H$ under
the action of any homomorphism of $\Psi $. The residual
$p$@-finiteness of group $G$ results that the centralizer $C$($H$)
of subgroup $H$ of group $A$ is a $p$@-separable subgroup, and so
$f \in C(H)$.

Obviously, if element $f$ has an infinite order modulo subgroup
$H$ (i.~e. $f^{n} \notin H$ for any natural $n$), it is
sufficient to put $f_{n} = f^{n}$. Therefore the order of $f$
modulo $H$ is considered to be finite and equal to $q$. We'll show
that $q$ isn't a $p$@-number.

Let $f^{q} = h^{m}$. Since, by the condition, subgroup $C(H)$
doesn't contain elements of finite order, $(m, q) = 1$. From the
order hand, there exists a homomorphism $\psi \in \Psi $, mapping
$h$ to a non-identity element. By virtue of the choice of element
$f$ one can find such a number $x$ that $f\psi = h^{x}\psi $.
Then $qx \equiv m~(mod~\vert h\psi \vert )$, and, the order of the
element $h\psi $ being a non-unit $p$@-number, the property
$p\vert q$ would imply $p\vert m$ and $(m, q) \ne 1$.

Thus, $q$ isn't a $p$@-number, and we can put $f_{n} = f^{n}$ again.

Let now $b \in B\smallsetminus K$ be an arbitrary element. Suppose
that $b^{- 1}Kb \cap K \ne 1$, and $b^{ - 1}k^{n}b \in b^{- 1}Kb
\cap K$ for some $n>0$. Since subgroup $K$ is separable by the
subgroups of family $\Omega _{B}^{p} $, it is $p^{\prime}
$@-isolated in group $B$, hence $n$ may be considered as a
$p$@-number. Putting $g = [b^{- 1}f_{n}b, f_{n}]$, where element
$f_{n}$ is defined above, we get $g \ne 1$ and at the same time
$g\psi = 1$ for any homomorphism $\psi \in \Psi $. This
contradicts the residual $p$@-finiteness of group $G$.

Thus, $b^{- 1}Kb \cap K = 1$ for every element $b \in
B\smallsetminus K$. It follows, in particular, that an arbitrary
non-unit element of subgroup $K$ doesn't commutate with any
element of group $G$ having a reduced form of a syllable length
greater than 1.

Let now a $p^{\prime} $@-isolated cyclic subgroup $\langle
u\rangle $ of group $G$ be not $p$@-separable in $G$, and let $v
\in G$ be such an element that $v \notin \langle u\rangle $, but
$v\psi \in \langle u\psi \rangle $ for every homomorphism $\psi
\in \Psi $. Then, as it was noted above, $[u, v] = 1$.

Suppose, at first, that element $u$ belongs to some subgroup $C$
conjugated with $A$ or with $B$. Then $v$ comes to be an element
of the same subgroup $C$: it follows from general considerations
(see, e.~g., [11, theorem~4.5]) if $u$ isn't contained in a
subgroup conjugated with $K$, and from proved before otherwise.
Hence, the subgroup $\langle u\rangle $ conjugates with a subgroup
of family $\Lambda _{A}^{p} \cup \Lambda _{B}^{p} $, which
coincides with $\Delta _{A}^{p} \cup \Delta _{B}^{p} $.

The case, when element $v$ belongs to a subgroup conjugated with
$A$ or with $B$, is considered similarly.

Let, at last, neither $u$, nor $v$ be contained in such a
subgroup. Then $u = g^{- 1}k^{m}gw^{s}$, $v = g^{ -
1}k^{n}gw^{t}$, where $g$,~$w \in G$ and $[g^{- 1}k^{m}g, w] =
[g^{- 1}k^{n}g, w] = 1$ [ibid.]. We'll show that it is impossible.

It follows from $[g^{- 1}k^{m}g, w] = [g^{- 1}k^{n}g, w] = 1$ that
$[k^{m}, gwg^{- 1}] = [k^{n}, gwg^{- 1}] = 1$, and so either $m =
n = 0$, or $gwg^{- 1} \in A \cup B$. The second case just gives a
contradiction, since elements $u = g^{ - 1}k^{m}(gwg^{ - 1})^{s}g$
and $v = g^{- 1}k^{n}(gwg^{- 1})^{t}g$ turn out in the subgroup
conjugated with $A$ or with $B$ by element $g$.

Thus, $u = w^{s}$ and $v = w^{t}$, and $s$ is a $p$@-number. As it
was noted above, the residual $p$@-finiteness of group $G$ gives
the existence of it's normal subgroup $N$, say, of finite
$p$@-index, which intersects with the cyclic subgroup generated by
element $w$ at the subgroup $\langle w^{s}\rangle $. Since $v
\notin \langle u\rangle $, $vN \ne 1$ in the group $G/N$, i.~e.
$vN \notin \langle uN\rangle $. We get a contradiction with the
choice of element $v$.

\bigskip

{\bf References}

\bigskip

\noindent\item{1.} {\it Mal'cev~A.~I.} On homomorphisms onto finite groups, {\it
Transl.\ Amer.\ Math.\ Soc.}, {\bf 119} (1983), 67--79; translation from
{\it Ivanov.\ Gos.\ Ped.\ Inst.\ Ucen.\ Zap.}, {\bf 18} (1958), 49--60
(Russian).

\noindent\item{2.} {\it Kim~G.} Cyclic subgroup separability of generalized free
products, {\it Canad.\ Math.\ Bull.}, (3) {\bf 36} (1993), 296--302.

\noindent\item{3.} {\it Baumslag~G.} On the residual finiteness of generalized
free products of nilpotent groups, {\it Trans.\ Amer.\ Math.\ Soc.},
{\bf 106} (1963), 193--209.

\noindent\item{4.} {\it Loginova~E.~D.} Residual finiteness of the free product of
two groups with commuting subgroups, {\it Sib.\ Math.~J.}, {\bf 40}
(1999), 341--350; translation from {\it Sib.\ Mat.\ Zh.} {\bf 40} (1999),
395--407 (Russian).

\noindent\item{5.} {\it Kim~G., Tang~C.~Y.} On generalized free products of
residually finite $p$@-groups, {\it J.~Algebra}, {\bf 201} (1998),
317--327.

\noindent\item{6.} {\it Azarov~D.~N.} On the residual nilpotence of free products
of free groups with cyclic amalgamation, {\it Math.\ Notes}, {\bf 64}
(1998), 3--7; translation from {\it Mat.\ Zametki} {\bf 64} (1998), 3--8
(Russian).

\noindent\item{7.} {\it Sokolov~E.~V.} The finite separability of cyclic subgroups
of some generalized free products of groups, {\it Vestnik Mol.\ Uchen.\ IvGU},
{\bf 2} (2002), 7--10 (Russian).

\noindent\item{8.} {\it Sokolov~E.~V.} On the residual $p$@-finiteness of some
free products with amalgamated subgroup, {\it Chebishevskii Sbornik},
(1) {\bf 3} (2002), 97--102 (Russian).

\noindent\item{9.} {\it Higman~G.} Amalgams of $p$@-groups, {\it J.~Algebra}, {\bf 1}
(1964), 301--305.

\noindent\item{10.} {\it Neumann~H.} Generalized free products with amalgamated
subgroups II, {\it Am.\ J.\ Math.}, {\bf 31} (1949), 491--540.

\noindent\item{11.} {\it Magnus~W., Karras~A., Solitar~D.} Combinatorial group
theory. Interscience Publishers, 1966; russian transl.: Moscow, Nauka, 1974.

\end